\title{ The largest eigenvalue of sparse random graphs
}
\author{Michael Krivelevich
\thanks{Department of Mathematics,
        Raymond and Beverly Sackler Faculty of
        Exact Sciences, Tel Aviv University, Tel Aviv 69978,
        Israel. Email: krivelev@math.tau.ac.il.
        Supported by a USA-Israeli BSF grant and by a Bergmann Memorial
        Award.
        }
\and
Benny Sudakov
\thanks{Department of Mathematics,
Princeton University, Princeton, NJ 08544, USA
and Institute for Advanced Study, Princeton, NJ 08540,
USA. Email address: bsudakov@math.princeton.edu.
Research supported in part by NSF grants DMS-0106589, CCR-9987845
and by the State of New Jersey.
}        }
\date{}
\documentclass [11pt]{article}
\usepackage{amsfonts}
\usepackage{amsmath}
\usepackage{latexsym}
\newtheorem{theo}{Theorem}[section]
\newtheorem{prop}[theo]{Proposition}
\newtheorem{lemma}[theo]{Lemma}
\newtheorem{coro}[theo]{Corollary}

\begin{document}
\maketitle

\begin{abstract}
We prove that for all values of the edge probability $p(n)$ 
the largest eigenvalue of a random graph $G(n,p)$ satisfies almost
surely:
$\lambda_1(G)=(1+o(1))\max\{\sqrt{\Delta},np\}$, where $\Delta$ is a 
maximal
degree of $G$, and the $o(1)$ term tends to zero as 
$\max\{\sqrt{\Delta},np\}$
tends to infinity. 
\end{abstract}

\section{Introduction}
Let $G=(V,E)$ be a graph with vertex set $V(G)=\{1,\ldots,n\}$. The
{\em adjacency matrix} of $G$, denoted by $A=A(G)$, is an $n$-by-$n$
$0,1$-matrix whose entry $A_{ij}$ is one if $(i,j)\in E(G)$, and is zero
otherwise. It is immediate that $A(G)$ is a real symmetric matrix. We
thus denote by $\lambda_1\ge\lambda_2\ge\ldots \lambda_n$ the
eigenvalues of $A$ which are usually called also the eigenvalues of the
graph $G$ itself. A family $\{\lambda_1,\ldots,\lambda_n\}$ is called
the {\em spectrum} of $G$.  

Spectral techniques play an increasingly important role in modern Graph
Theory. A serious effort has been invested in establishing connections
between spectral characteristics of a graph and its other parameters. An
interested reader can consult monographs \cite{CveDooSac95},
\cite{Chu97} for a detailed account of known results. The ability to
compute graph eigenvalues efficiently (both from theoretical and
practical points of view), combined with results from spectral graph
theory, has provided a basis for quite a few graph algorithms. A survey
of applications of spectral techniques in Algorithmic Graph Theory by
Alon can be found in \cite{Alo98}.

In this paper we study eigenvalues of random graphs. A {\em random
graph} $G(n,p)$ is a discrete probability space composed of all labeled
graphs on vertices $\{1,\ldots,n\}$, where each edge $(i,j)$,
$1\le i<j\le n$, appears randomly and independently with probability 
$p=p(n)$.
Sometimes with some abuse of notation we will refer to a random graph
$G(n,p)$ as a graph on $n$ vertices generated according to the
distribution $G(n,p)$ described above.
Usually asymptotic properties of random graphs are of interest. We say
that a graph property ${\cal A}$ holds {\em almost surely}, or a.s. for
brevity, in $G(n,p)$ if the probability that $G(n,p)$ has ${\cal A}$
tends to one as the number of vertices $n$ tends to infinity. Necessary
background information on random graphs can be found at \cite{B},
\cite{JanLucRuc00}. It is
important to observe that the adjacency matrix of a random graph
$G(n,p)$ can be viewed as  a random symmetric matrix, whose diagonal
entries are zeroes and whose entries above the diagonal are i.i.d.
random variables, each taking value 1 with probability $p$ and value 0
with probability $1-p$. This allows to bridge between random graphs and
extensively developed theory of random real symmetric matrices and
their spectra (see, e.g. \cite{Meh91}). 

The subject of this paper is asymptotic behavior of the largest
eigenvalue $\lambda_1(G(n,p))$ of random graphs. (Notice that due to the
Perron-Frobenius Theorem, for every graph $G$ on $n$ vertices,
$\lambda_1(G)\ge|\lambda_i(G)|$ for all $i=2,\ldots,n$.
Thus $\lambda_1(G)$ is equal to the spectral norm or the spectral radius
of $A(G)$. It is easy to observe that for every graph $G=(V,E)$ its
largest eigenvalue $\lambda_1(G)$ is always squeezed between the average
degree of $G$, $\bar{d}=\sum_{v\in V}d_G(v)/|V|$ and its maximal degree
$\Delta(G)=\max_{v\in V}d_G(v)$. As for all $p(n)\gg \log n$ the last
two quantities are both asymptotically equal to $np$, it follows that in
this range of edge probabilities a.s. $\lambda_1(G(n,p))=(1+o(1))np$. In
fact, much more is known for large enough values of $p(n)$, F\"uredi and
Koml\'os proved in \cite{FurKom81} that for a constant $p$,
$\lambda_1(G(n,p))$ has asymptotically a normal distribution with
expectation $(n-1)p+(1-p)$ and variance $2p(1-p)$.

In contrast, not much appears to be known for the case of sparse random
graphs, i.e. when $p(n)=O(\log n)$. Khorunzhy and Vengerovsky
\cite{KhoVen00} and Khorunzhy \cite{Kho01b} consider mainly the case
$p(n)=1/n$ and show that in this case the spectral norm of $A(G(n,p))$
a.s. tends to infinity with $n$. Moreover, it is stated in
\cite{KhoVen00} that the mathematical expectation of the number of
eigenvalues that go to infinity is of order $\Theta(n)$.      

Here we are to find the asymptotic value of the largest eigenvalue of
sparse random graphs. To grasp better the result, observe that if
$\Delta$ denotes a maximal degree of a graph $G$, then $G$ contains a
star $S_\Delta$ and therefore $\lambda_1(G)\ge
\lambda_1(S_\Delta)=\sqrt{\Delta}$. Also, as mentioned above
$\lambda_1(G)$ is at least as large as an average degree of $G$. As for
all values of $p(n)\gg 1/n^2$, a.s. $|E(G(n,p))|=(1+o(1))(n^2p/2)$, we
get that a.s. $\lambda_1(G(n,p))\ge (1+o(1))np$. Combining the above
lower bounds together, we get that a.s. $\lambda_1(G(n,p))\ge
(1+o(1))\max\{\sqrt{\Delta}, np\}$. As it turns out this lower bound can
be matched by an upper bound of the same asymptotic value, as given by
the following theorem:

\begin{theo} \label{t11}
Let $G=G(n,p)$ be a random graph and let $\Delta$ be a maximum degree of
$G$. Then almost surely the largest eigenvalue of the adjacency matrix
of $G$ satisfies
$$
 \lambda_1(G)=(1+o(1))\max\big\{\sqrt{\Delta}, np\big\},
$$
where $o(1)$ tends to zero as $\max\{\sqrt{\Delta}, np\}$ tends to
infinity.
\end{theo}

As the asymptotic value of the maximal degree of $G(n,p)$ is known for
all values of $p(n)$ (see Lemma \ref{l21} below), the above theorem
enables to estimate the asymptotic value of $\lambda_1(G(n,p))$ for all
relevant values of $p$. In particular, for the case $p=c/n$ we get:

\begin{coro}\label{cor1}
For any constant $c>0$, a.s. 
$\lambda_1(G(n,c/n))=(1+o(1))\sqrt{\frac{\ln n}{\ln\ln n}}$.
\end{coro}

The rest of the paper is organized as follows. In the next section we
gather necessary technical information about random graphs, used later
in the proof of the main result. The main theorem, Theorem \ref{t11},
is proven
in Section 3. Section 4, the last section of the paper, is devoted to
concluding remarks and discussion of related open problems.

Throughout the paper we will omit systematically floor and ceiling
signs for the sake of clarity of presentation. All logarithms are
natural.  We will frequently use the inequality
 ${a\choose b}\le \left(\frac{ea}{b}\right)^b$.

\section{Some properties of sparse random graphs}
In this section we show some properties of sparse random graphs
which we will use later to prove Theorem \ref{t11}. First we need the
following definition. Let $G(n,p)$ be a random graph. Denote
$$
\Delta_p=\max\Big\{k: n{{n-1}\choose k}p^k(1-p)^{n-k}\ge 1\Big\}\ .
$$
In words, $\Delta_p$ is the maximal $k$ for which the expectation of
the number of vertices of degree $k$ in $G(n,p)$ is still at least one. 

\begin{lemma}
\label{l21}
Let $G=G(n,p)$ be a random graph. Then

\noindent
$(i)$\, The maximum degree of $G$ almost surely satisfies
$\Delta(G)=(1+o(1))\Delta_p$.

\noindent
$(ii)$\, If $np \rightarrow 0$ then almost surely $G$ is a forest.

\noindent
$(iii)$\, If $p \leq e^{-(\log \log n)^2}/n$, then
almost surely all connected
components of $G$ are of size at most $(1+o(1))\Delta_p$.

\noindent
$(iv)$\,  If $p \leq \log^{1/2} n/n$, then almost surely every vertex of
$G$ is contained in at most one cycle of length $\leq 4$.
\end{lemma}

\noindent
{\bf Proof.}\, Parts (i) and (ii) are well known and can be found,
e.g., in the monograph of Bollob\'as \cite{B}. To show (iii) it is enough to
bound from above the expectation of the number $Y$ of trees on
$t=(1+1/\log \log n)\Delta_p +1$ vertices, contained in $G(n,p)$ as
subgraphs. Obviously this expectation is equal to
$$EY={n \choose t}t^{t-2}p^{t-1} \leq
\frac{n^t}{t!}t^{t-2}p^{t-1} \leq \frac{n^t}{(t/e)^t}t^{t-2}p^{t-1}=
\frac{en}{t^2}\big(enp\big)^{t-1}=\frac{e}{t^2}\Big(n
\big(enp\big)^{\Delta_p}\Big)
\big(enp\big)^{t-1-\Delta_p}.$$
On the other hand, by the definition of $\Delta_p$, we have that
$n\big(enp\big)^{\Delta_p}=O({\Delta_p}^{\Delta_p+1})$ and
${\Delta_p}=o(\log n)$. Therefore, using that $p \leq e^{-(\log \log
n)^2}/n$  and $t>\Delta_p$, we conclude
$$EY \leq O\left(\frac{e}{t^2}
{\Delta_p}^{\Delta_p+1}\big(enp\big)^{\Delta_p/\log \log n}\right)
\leq O\left( \Big(\frac{e\Delta_p}{\log n}\Big)^{\Delta_p}\right)=o(1).$$
Now (iii) follows from Markov's inequality.
Finally, the expected number of pairs of intersecting cycles of length $s,
t \leq 4$ in the graph $G$ is obviously at most
$O(n^sn^{t-1}p^{s+t}) \leq O(\log^4 n/n)=o(1)$. This, by Markov's
inequality, implies (iv).
$\quad \Box$

\bigskip

Next we show that the set of vertices of relatively high degree in
$G(n,p)$ spans a graph with small maximum degree and with
no cycles. More precisely, the following stronger
statement is true.

\begin{lemma}
\label{l22}
Let $p \geq e^{-(\log \log n)^2}/n$ and let $X$ be the set of vertices of
random graph
$G=G(n,p)$ with
degree larger than $np(1+1/\log \log n)+\Delta_p^{1/3}$.
Then

\noindent
$(i)$\, Almost surely every cycle of $G$ of length $k$ intersects
$X$ in less than $k/2$ vertices.

\noindent
$(ii)$\, Almost surely every vertex in $G$ has less than
$\Delta_p^{7/8}$ neighbors in $X$.
\end{lemma}

\noindent
{\bf Proof.}\, First we consider the case when
$e^{-(\log \log n)^2}/n \leq p \leq \log^{1/4} n/n$. In this case, by
definition, $\Delta_p=\Omega(\log n/(\log \log n)^2)$ and
$np \leq \log^{1/4} n$. To prove the lemma we first estimate the
probability that all the vertices of  a fixed set $T$ of size
$|T|=t$ have
degrees at least $\log^{1/3} n/\log \log n < \Delta_p^{1/3}$. It is easy
to see that for such a set $T$, either there are at least
$(\log^{1/3} n/\log \log n)t/3$ edges in the cut $(T,V(G)-T)$,
or the set $T$ spans at least $(\log^{1/3} n/\log \log n)t/3$ edges of
$G$. Since the number of edges in the cut
$(T,V(G)-T)$ is a binomially distributed random variable with parameters
$t(n-t)$ and $p$,
we can bound the probability of the first event by

\begin{eqnarray*}
{t(n-t) \choose \frac{\log^{1/3} n}{3\log \log n}t}
p^{\frac{\log^{1/3} n}{3\log \log n}t}&\leq&
\left(\frac{3e(n-t)p\log \log n}{\log^{1/3} n}\right)^{\frac{\log^{1/3}
n}{3\log \log n}t}\leq
\left(\frac{3e \log^{1/4} n \log \log n}{\log^{1/3}
n}\right)^{\frac{\log^{1/3} n}{3\log \log n}t}\\
&\leq&
e^{-\Omega(t\log^{1/3} n)}.
\end{eqnarray*}
Also, the number of edges spanned
by $T$ is a binomially distributed random variable with parameters
$t(t-1)/2$ and $p$. We can thus bound
the probability of the second event similarly by
\begin{eqnarray*}
{\frac{t(t-1)}{2} \choose \frac{\log^{1/3} n}{3\log \log n}t}
p^{\frac{\log^{1/3} n}{3\log \log n}t}&\leq&
\left(\frac{3e(t-1)p\log \log n}{2\log^{1/3} n}\right)^{\frac{\log^{1/3}
n}{3\log \log n}t}\leq
\left(\frac{3e \log^{1/4} n \log \log n}{2\log^{1/3}
n}\right)^{\frac{\log^{1/3} n}{3\log \log n}t}\\
&\leq&
e^{-\Omega(t\log^{1/3} n)}.
\end{eqnarray*}
Therefore, the probability that all the vertices in the given
set of size $t$ have degree at least $\Delta_p^{1/3}$ is at most
$e^{-\Omega(t\log^{1/3} n)}$. Essentially repeating the above argument
shows that conditioning on the presence of any specific set of at most
$2t$ edges in $G$ leaves the latter probability still at most
$e^{-\Omega(t\log^{1/3}n)}$.

Using this bound we can easily estimate the probability that there
exists a cycle
of length $k$ with at least $k/2$ vertices inside the set $X$.
Clearly this probability is at most
$$\sum_{k \geq3} n^kp^k {k \choose {\lceil k/2\rceil}} e^{-\Omega((k/2)\log^{1/3} n)}
\leq \sum_{k \geq3}
\left(2np e^{-\Omega(\log^{1/3} n)}\right)^k
\leq \sum_{k \geq3}
\left(2(\log^{1/4} n) e^{-\Omega(\log^{1/3} n)}\right)^k=o(1).$$
(First choose  $k$ vertices of a cycle and fix their order, then require
that the $k$ edges of the cycle are present in $G(n,p)$, then choose a
set $T$ of the cycle vertices of cardinality $|T|=t=\lceil k/2\rceil$,
and then require all vertices of $T$ to belong to $X$, conditioning on
the presence of the cycle edges in $G(n,p)$.)
This implies claim (i) of the lemma.
Similarly, the probability that there exists a vertex with at least
$\Delta_p^{7/8}$ neighbors in $X$ is at most
\begin{eqnarray*}
n {n \choose \Delta_p^{7/8}} p^{\Delta_p^{7/8}}
e^{-\Omega(\Delta_p^{7/8}\log^{1/3} n)} \hspace{-0.1cm}&\leq&
\hspace{-0.1cm}
n\left(np e^{-\Omega(\log^{1/3} n)}\right)^{\Delta_p^{7/8}} \leq
n \left((\log^{1/4} n) e^{-\Omega(\log^{1/3} n)}\right)^
{\Omega\big((\frac{\log n}{(\log \log n)^2})^{7/8}\big)}\\
 \hspace{-0.1cm} &\leq& \hspace{-0.1cm}
n e^{-\Omega(\log^{13/12} n)}=o(1).
\end{eqnarray*}
This completes the proof of the lemma for
$e^{-(\log \log n)^2}/n \leq p \leq \log^{1/4} n/n$.

Next we consider the case when $p \geq \log^{1/4} n/n$. We again start
by
estimating the probability that
that all the vertices of a fixed set $T$ of size $t \leq n/2$ have
degree at least $np(1+1/\log \log n)$. Similarly as before,
for such a set $T$, there are at least
$t(n-t)p+tnp/(3\log \log n)$ edges in the cut $(T,V(G)-T)$, or the set
$T$ spans at least $t(t-1)p/2+tnp/(3\log \log n)$ edges.
By the standard estimates for Binomial distributions
(see, e.g., \cite{AS}, Appendix A) it follows that the probability of
the first event is at most $e^{-\Omega(tnp/(\log \log n)^2)}$. The same
estimates can be used to show that if $n/(6\log \log n) \leq t \leq n/2$
then the probability of the second event is also bounded by
$e^{-\Omega(tnp/(\log \log n)^2)}$. On the other hand, if
$t \leq n/(6\log \log n)$, then this probability can be bounded directly by
$${\frac{t(t-1)}{2} \choose \frac{tnp}{3\log \log n}}p^{\frac{tnp}{3\log
\log n}}\leq \left(\frac{3e(t-1)p\log \log n}{2np}
\right)^{\frac{tnp}{3\log \log n}} \leq
\left(\frac{e}{4}\right)^{\frac{tnp}{3\log \log n}} \leq
e^{-\Omega(tnp/(\log \log n)^2)}.$$
Therefore, the probability that all degrees of vertices in the given
set of size $t$ are at least $np(1+1/\log \log n)$ is at most
$e^{-\Omega(tnp/(\log \log n)^2)}$. Again, conditioning on the presence of
any specific set of at most $2t$ edges does not change the order of the
exponent in the above estimate.

Using this bound together with the fact that
$np \geq \log^{1/4} n$ we can estimate probability that
there exists a cycle
of length $k$ with at least $k/2$ vertices inside set $X$.
Clearly this probability is at most
$$\sum_{k \geq3} n^kp^k {k \choose {\lceil k/2\rceil}}
e^{-\Omega((k/2)np/(\log \log n)^2)}
\leq \sum_{k \geq3}
\left(2np e^{-\Omega(np/(\log \log n)^2)}\right)^k
\leq \sum_{k \geq3} e^{-\Omega\big(k\frac{\log^{1/4} n}{(\log \log
n)^2}\big)}=o(1).$$
This implies claim (i).
Similarly, the probability that there exists a vertex with at least
$\Delta_p^{7/8}$ neighbors in $X$ is at most
\begin{eqnarray*}
n {n \choose \Delta_p^{7/8}} p^{\Delta_p^{7/8}}
e^{-\Omega(\Delta_p^{7/8}np/(\log \log n)^2)} &\leq&
n\left(np e^{-\Omega(np/(\log \log n)^2)}\right)^{\Delta_p^{7/8}} \leq
n e^{-\Omega\big(\frac{\log^{1/4} n}{(\log \log n)^2}
(\frac{\log n}{(\log \log n)^2})^{7/8}\big)}\\
&\leq&
n e^{-\Omega(\log^{17/16} n)}=o(1).
\end{eqnarray*}
This implies claim (ii) and completes the proof of the lemma.
$\quad \Box$

\bigskip

Finally we need one additional lemma.
\begin{lemma}
\label{l23}
Let $G=G(n,p)$ be a random graph with $e^{-(\log \log n)^2}/n \leq p
\leq \log^{1/2} n/n$.
Then a.s. $G$ contains no vertex which has
at least $\Delta_p^{1/3}$
other vertices of $G$ with degree $\geq \Delta_p^{3/4}$
within distance at most two.
\end{lemma}

\noindent
{\bf Proof.}\, Let $v$ be a vertex of $G(n,p)$ and let $u_i, i=1, \ldots,
\Delta_p^{1/3}$ be the vertices with degree at least $\Delta_p^{3/4}$
which are within distance at most two from $v$. Let $T$ be the set of
vertices of the smallest connected  subgraph of $G$ which contain $v$
together with all the vertices $u_i$. Since the shortest path form $v$ to
$u_i$ may contain only one vertex distinct from $v$ and $u_i$, then it is
easy to see that the size
of $T$ satisfies $\Delta_p^{1/3}+1 \leq |T|=t \leq 2\Delta_p^{1/3}+1$. In
addition each $u_i$ has at least $\Delta_p^{3/4}-t \geq
\frac{1}{2}\Delta_p^{3/4}$
neighbors outside set $T$. Therefore there are at least
$\frac{1}{2} \Delta_p^{3/4} \cdot
\Delta_p^{1/3}=\frac{1}{2} \Delta_p^{13/12}$ edges of $G$ between
$T$ and $V(G)-T$. Since the number of edges in the cut
$(T,V(G)-T)$ is a binomially distributed random variable with parameters
$t(n-t)$ and $p$
we can bound the probability of this event for a fixed set $T$ of size
$|T|=t$ by
$${t(n-t) \choose
{\frac{1}{2}\Delta_p^{13/12}}}p^{\frac{1}{2}\Delta_p^{13/12}} \leq
\left(\frac{2et(n-t)p}{\Delta_p^{13/12}}\right)^{\frac{1}{2}\Delta_p^{13/12}}
\leq \left(\frac{5e \log^{1/2}
n}{\Delta_p^{3/4}}\right)^{\frac{1}{2}\Delta_p^{13/12}} \leq
e^{-\log^{25/24}n}.$$
Here we used that for $p \geq e^{-(\log \log n)^2}/n$, by definition,
$\Delta_p \geq \Omega(\log n/(\log \log n)^2)$  and the facts that $np \leq
\log^{1/2} n$ and $t \leq 2\Delta_p^{1/3}+1$.

As we explained in the previous paragraph, the probability that there
exists a vertex that  violates the assertion of the lemma is bounded by
the probability that  there exists a connected subgraph on $|T|=t \leq
2\Delta_p^{1/3}+1$  vertices such that the number of edges in the cut
$(T,V(G)-T)$ is at least $\frac{1}{2}\Delta_p^{13/12}$. Using that
for $p \leq \log^{1/2} n/n$, by definition,
$\Delta_p=o(\log n)$, we can bound this probability by
\begin{eqnarray*}
\sum_{t\leq 2\Delta_p^{1/3}+1} {n \choose t} t^{t-2} p^{t-1}
e^{-\log^{25/24}n}\hspace{-0.1cm} &\leq& \hspace{-0.2cm}
\sum_{t\leq 2\Delta_p^{1/3}+1}\frac{en}{t^2}
\big(enp\big)^{t-1}e^{-\log^{25/24}n} \leq
3\Delta_p^{1/3}n\big(enp\big)^{2\Delta_p^{1/3}}e^{-\log^{25/24}n}
\\
\hspace{-0.1cm}&\leq&
n \log^{1/3} n \big(e\log^{1/2} n\big)^{\log^{1/3}n}
e^{-\log^{25/24}n}=o(1).
\end{eqnarray*}
This completes the proof. $\quad \Box$

\section{The proof of main result}
In this section we prove Theorem \ref{t11}. We start by stating
some simple properties of the largest eigenvalue of a graph,
that we will need later.

\begin{prop}
\label{p31}
Let $G$ be a graph on $n$ vertices and $m$ edges and with maximum degree
$\Delta$. Let $\lambda_1(G)$ be the largest eigenvalue of the
adjacency matrix of $G$. Then is has the following properties.

\noindent
$(i)$\, $\max\big(\sqrt{\Delta}, \frac{2m}{n}\big) \leq
\lambda_1(G) \leq \Delta$.

\noindent
$(ii)$\, If $E(G)= \cup_i E(G_i)$ then $\lambda_1(G) \leq \sum_i
\lambda_1(G_i)$. If in addition graphs
$G_i$ are vertex disjoint, then $\lambda_1(G)=\max_i \lambda_1(G_i)$.

\noindent
$(iii)$\, If $G$ is a forest then $\lambda_1(G) \leq \min
\big(2\sqrt{\Delta-1}, \sqrt{n-1}\big)$. In particular if $G$ is a star
on $\Delta+1$ vertices then $\lambda_1(G)=\sqrt{\Delta}$.

\noindent
$(iv)$\, If $G$ is a bipartite graph such that degrees on both sides of
bipartition are
bounded by $\Delta_1$ and $\Delta_2$ respectively, then
$\lambda_1(G) \leq \sqrt{\Delta_1 \Delta_2}$.
\end{prop}

\noindent
{\bf Proof.}\, Most of these easy statements can be found in Chapter 11
of the book of Lov\'asz \cite{L}. Here we sketch the proof of few
remaining ones for the sake of completeness.

\noindent
(iii)\, Let $A$ be the adjacency matrix of $G$ and let
$\lambda_1, \lambda_2, \ldots, \lambda_n$ be its eigenvalues.
Since $G$ is a forest on $n$
vertices it is easy to see that the trace of $A^2$ satisfies
$$tr(A^2)=\sum_i \lambda_i^2=\sum_v d_v\leq 2(n-1).$$
On the other hand $\lambda_1=-\lambda_n$ because $G$ is bipartite.
Therefore we can conclude that $2\lambda_1^2\leq 2(n-1)$ and
hence $\lambda_1 \leq \sqrt{n-1}$. For the proof of the rest of the
statement (iii) see, e.g., \cite{L}.

\noindent
(iv)\, Let $A$ be the adjacency matrix of $G$. Then by definition it is
easy to see that $A^2$ is the adjacency matrix of a multigraph with maximum
degree $\Delta_1 \Delta_2$. Therefore by (i) we have that
$\lambda_1(A^2)=\lambda_1^2(G) \leq \Delta_1 \Delta_2$ and
hence $\lambda_1 \leq \sqrt{\Delta_1 \Delta_2}$. $\quad \Box$

\bigskip

Having finished all the necessary preparations, we are now ready to
complete the proof of our main theorem.

\noindent
{\bf Proof of Theorem \ref{t11}.}\,
We start with the easy case when the random graph is very sparse.
If $p \leq e^{-(\log \log n)^2}/n$, then by Lemma
\ref{l21} a.s. $G=G(n,p)$ is a disjoint union of trees of size at most
$(1+o(1))\Delta_p$.
Therefore by
claims (ii) and (iii) of
Proposition \ref{p31} we have that $\lambda_1(G) \leq
(1+o(1))\sqrt{\Delta_p}$. On the other hand, By Lemma \ref{l21} the
maximum degree of $G$ is almost surely at least
$(1+o(1))\Delta_p$ and thus claim (i)
of Proposition \ref{p31} implies that
$\lambda_1(G) \geq (1+o(1))\sqrt{\Delta_p}$. Since
the value of the edge probability satisfies $np=o(1)< \Delta_p$,
 we obtain that
$\lambda_1(G)=(1+o(1))\sqrt{\Delta_p}=(1+o(1))
\max\big(\sqrt{\Delta(G)}, np\big)$.

Another relatively simple case is when $p \geq \log^{1/2} n/n$. Then by
definition it is easy to check that $\Delta_p=o\big((np)^2\big)$ and hence
it is enough to prove that
$\lambda_1(G)=(1+o(1))\max\big(\sqrt{\Delta_p}, np\big)=
(1+o(1))np$. To get a lower bound on the largest eigenvalue note that
the standard Chernoff estimates for the binomial distributions
(see, e.g., \cite{AS} , Appendix A) imply that the
number of edges in $G(n,p)$ is a.s. $(1+o(1))n^2p/2$. Therefore by
claim (i) of Proposition \ref{p31},  the largest eigenvalue of
$G(n,p)$ is almost surely  at least $(1+o(1))n^2p/n=(1+o(1))np$.

To get an upper bound,  denote by $X$ the set of vertices of
random graph $G=G(n,p)$ with
degree larger than $np(1+1/\log \log n)+\Delta_p^{1/3}$.
Let $G_1$ be a subgraph of $G$ induced by the set $X$, let $G_2$ be
a subgraph of $G$ induced by the set $V(G)-X$ and finally let
$G_3$ be a bipartite subgraph of $G$ containing all the edges between $X$
and $V(G)-X$. By definition $G=\cup_i G_i$ and thus by claim (ii) of
Proposition \ref{p31} we obtain that
$\lambda_1(G) \leq \sum_{i=1}^3 \lambda_1(G_i)$.
Since the maximum degree of graph $G_2$ is $np(1+1/\log \log
n)+\Delta_p^{1/3}=(1+o(1))np$, then by claim (i) of
Proposition \ref{p31} it follows that $\lambda_1(G_2) \leq (1+o(1))np$.
Also note that by our construction, any cycle in the graphs $G_1$ or $G_3$
should have at least half of its vertices in the set $X$. Therefore
from Lemma \ref{l22} we get that almost surely $G_1$ and $G_3$ contains no
cycles. In addition, by Lemma \ref{l21}, the maximum degree of these two
forests is
bounded by $(1+o(1))\Delta_p$. Then using claim (iii) of Proposition
\ref{p31} we
obtain that $\lambda_1(G_i) \leq (2+o(1))\sqrt{\Delta_p},~i=1,3$.
This implies that
$$\lambda_1(G) \leq \lambda_1(G_1)+\lambda_1(G_2)+\lambda_1(G_3)
\le (1+o(1))np+(4+o(1))\sqrt{\Delta_p}=(1+o(1))np.$$

Finally we treat the remaining case when
$e^{-(\log \log n)^2}/n \leq p \leq  \log^{1/2} n/n$.
Similarly as before we have that a.s. the maximum degree of $G=G(n,p)$ is
$(1+o(1))\Delta_p$ and the total number of edges in $G$ is
$(1+o(1))n^2p/2$. Therefore claim (i) of Proposition \ref{p31} implies
that
$$\lambda_1(G) \geq (1+o(1))\max\big(\sqrt{\Delta_p}, n^2p/n\big)=
 (1+o(1))\max\big(\sqrt{\Delta(G)}, np\big).$$
To handle the upper bound on $\lambda_1$ we again use a partition of $G$
into smaller subgraphs, whose largest eigenvalue is easier to
estimate.

Denote by $X_1$ the set of vertices of $G$ with degree at least
$\Delta_p^{3/4}$ and by $X_2$ the set of vertices with degrees larger
than $np(1+1/\log \log n)+\Delta_p^{1/3}$ but less than $\Delta_p^{3/4}$.
Let $X=X_1 \cup X_2$ and let $Y_1$ contains all vertices of $V(G)-X$
with at least one neighbor in $X_1$. Finally let $Y_2$ be the set
$V(G)-X\cup Y_1$. Note that by definition there are no edges between
$X_1$ and $Y_2$.

 We consider the following subgraphs of $G$.
Let $G_1$ be the subgraph of $G$ induced by the set $X$. Then by
Lemma \ref{l22}, $G_1$ contains no cycles and has maximum degree
at most $\Delta_p^{7/8}$. Therefore by claim (iii) of
Proposition \ref{p31} we get that
$\lambda_1(G_1) \leq
2\sqrt{\Delta_p^{7/8}}=o\big(\sqrt{\Delta_p}\big)$.

Our second graph $G_2$ consists of all edges between $X_2$ and
$V(G)-X$. Note that by definition, any cycle in $G_2$  has exactly half
of its vertices  in $X_2 \subset X$. Thus by Lemma \ref{l22}, almost
surely $G_2$ is a forest. In addition, the maximum degree in $G_2$ is
bounded by the maximal possible degree of a vertex from the set
$V(G)-X_1$, which is $\Delta_p^{3/4}$. Using claim (iii) of
Proposition \ref{p31} we get that $\lambda_1(G_2) \leq
2\sqrt{\Delta_p^{3/4}}=o\big(\sqrt{\Delta_p}\big)$.

Next consider the graph $G_3$, induced by the set of vertices $Y_1$.
Let $v \in V(G)-X$ be a vertex with at least
$\Delta_p^{1/3}+1$ neighbors in $Y_1$.
Since by definition every neighbor of $v$ in $Y_1$ is also a neighbor
of some vertex in $X_1$ we obtain that there are at least
$\Delta_p^{1/3}+1$ paths of length two from $v$ to the set $X_1$.
On the other hand, by Lemma \ref{l21}, $v$ almost surely is contained
in at most one cycle of length 4. This implies that all but at most
one of the endpoints of these paths in $X_1$ are different. Therefore
vertex $v$ has at least $\Delta_p^{1/3}$ distinct vertices of
$X_1$ with in distance two. Now from Lemma \ref{l23}
it follows that a.s. there is no vertex with this property. Hence every
vertex  $v \in V(G)-X$ has almost surely
at most  $\Delta_p^{1/3}$ neighbors in $Y_1$. In particular the maximum
degree of $G_3$ is bounded by $\Delta_p^{1/3}$, which implies that
$\lambda_1(G_3) \leq \Delta_p^{1/3}=
o\big(\sqrt{\Delta_p}\big)$.

Let $G_4$ be the bipartite subgraph consisting of all the edges of $G$
between $Y_1$ and $Y_2$. By definition, the degree of every vertex in
$Y_1$ is at most $np(1+1/\log \log n)+\Delta_p^{1/3}$ and we already proved
in the previous paragraph that the degree of every vertex from $Y_2$ in
this graph is at most $\Delta_p^{1/3}$. Therefore, using claim
(iv) of Proposition \ref{p31} together with the facts that
$np \leq \log^{1/2} n$ and $\Delta_p=\Omega\big(\log n/(\log \log
n)^2\big)$ we obtain
$$\lambda_1(G_4) \leq \sqrt{\Delta_p^{1/3}\big(np(1+1/\log \log
n)+\Delta_p^{1/3}\big)} \leq \Delta_p^{1/3}+(1+o(1))\Delta_p^{1/6}
\sqrt{np}=o\big(\sqrt{\Delta_p}\big).$$

Finally we define $G_5$ to be the subgraph of $G$ induced by the set
$Y_2$, and $G_6$ to be a bipartite graph containing all the edges of
$G$ between $X_1$ and $Y_1$. Since there are no edges crossing from
$X_1$ to $Y_2$ it is easy to check that $E(G)=\bigcup_{i=1}^6 E(G_i)$.
Also since
the graphs $G_5$ and $G_6$ are vertex disjoint, then by claim (ii)
of Proposition \ref{p31} we obtain that
$\lambda_1(G_5 \cup G_6)=\max\big(\lambda_1(G_5), \lambda_1(G_6)\big)$
and
$$\lambda_1(G) \leq \lambda_1(G_1)+ \ldots + \lambda_1(G_4)+
\lambda_1(G_5 \cup G_6)=\max\big(\lambda_1(G_5), \lambda_1(G_6)\big)+
o\big(\sqrt{\Delta_p}\big).$$
By definition, the maximum degree of $G_5$ is bounded by
$(1+o(1))np+\Delta_p^{1/3}$, which implies that
$\lambda_1(G_5) \leq (1+o(1))np+\Delta_p^{1/3}$. Hence to finish the
proof it remains to bound $\lambda_1(G_6)$.

Consider the graph $G_6$.
Let $T$ be the set of vertices from $Y_1$ with degrees greater than one
in $G_6$ and let $u \in X_1$ be a vertex with at least
$\Delta_p^{1/3}+1$ neighbors in $T$.
By definition, every neighbor of $u$ in $T$ has also an additional
neighbor in $X_1$, which is distinct from $u$.
Therefore we obtain that there are
at least
$\Delta_p^{1/3}+1$ simple paths of length two from $u$ to the set $X_1$.
On the other hand, by Lemma \ref{l21}, $u$ almost surely is contained
in at most one cycle of length 4. This implies that all but at most
one of the endpoints of these paths in $X_1$ are different. Therefore
vertex $u$ has at least $\Delta_p^{1/3}$ distinct vertices of
$X_1$ within distance two. Now from Lemma \ref{l23}
it follows that a.s. there is no vertex with this property.
In addition it follows that every
vertex from $Y_1$ has degree at most $\Delta_p^{1/3}$ in $G_6$.
Let $H$ be the subgraph of $G_6$ containing all the edges from $X_1$ to
$T$. Then by the above discussion its maximum degree is bounded by
$\Delta_p^{1/3}$ and therefore $\lambda_1(H) \leq \Delta_p^{1/3}$.
On the other hand, since the degree of every vertex in $Y_1-T$ in $G_6$
is at most one and the graph is bipartite, we obtain that $G_6-H$ is a
union of vertex disjoint stars.
The size of each star is at most the maximum degree of $G$. Then by
claims (ii) and (iii) of Proposition \ref{p31} we get that
$$\lambda_1(G_6) \leq \lambda_1(H)+\lambda_1(G_6-H) \leq
\Delta_p^{1/3}+(1+o(1))\sqrt{\Delta_p}.$$
This implies the desired upper bound on $\lambda_1(G)$, since
\begin{eqnarray*}
\lambda_1(G) &\leq& \max\Big(\lambda_1(G_5), \lambda_1(G_6)\Big)+
o\Big(\sqrt{\Delta_p}\Big)=
\max\Big((1+o(1))np+\Delta_p^{1/3},
(1+o(1))\sqrt{\Delta_p}+\Delta_p^{1/3}\Big)\\
&=&
(1+o(1))\max\Big(\sqrt{\Delta_p}, np\Big)=
(1+o(1))\max\Big(\sqrt{\Delta(G)}, np\Big),
\end{eqnarray*}
and completes the proof of the theorem. $\quad \Box$

\section{Concluding remarks}
In this paper we have found the asymptotic value of the largest eigenvalue 
of a random
graph $G(n,p)$, or the spectral radius of the corresponding random real
symmetric matrix.

It would be quite interesting to obtain more accurate estimates on the
error term in the asymptotic estimate for $\lambda_1(G(n,p))$, given by
Theorem \ref{t11}. Notice that due to the recent concentration result of
the first author and Vu \cite{KriVu00}, the standard deviation of
$\lambda_1(G(n,p))$ can be asymptotically bounded by an absolute
constant, and this random variable is sharply concentrated. Our proof
methods do not allow us to locate the expectation of $\lambda_1$ with
such degree of precision. Neither we are able to obtain a limit
distribution of $\lambda_1$, as has been done by F\"uredi and Koml\'os
\cite{FurKom81}
for the case of a constant edge probability $p$. This is another
attractive open question.

One can also try to determine when the largest eigenvalue of 
a random graph has multiplicity one and then to 
understand a typical structure of the first
eigenvector of $G(n,p)$. While for the case $p\gg \log n/n$, where the
graph $G(n,p)$ becomes a.s. almost regular, the first eigenvector will
be a.s. almost collinear to the all-1 vector, the picture becomes more
complicated for smaller values of $p(n)$. Notice that for
$p(n)\ll \log n/n$ the graph $G(n,p)$ is a.s. disconnected, and
therefore the support of the first eigenvector will be at most as large
as the size of its largest connected component.

Consider the case $p=c/n$, for a constant $c>0$. Performing direct
calculations similar to those of Section 2 of the present paper, one 
can show that in this case $G(n,p)$ contains almost surely an unbounded
collection of vertices of degree $\Delta(G)(1-o(1))$ at distance at
least three from each other. Considering then the subgraph of $G$
 spanned by those vertices and their neighbors shows that a.s. $G(n,p)$
has an unbounded number of eigenvalues $\lambda_i=(1-o(1))\lambda_1$.

Another observation for the case $p=c/n$ is that according to Corollary
\ref{cor1} the first eigenvalue of $G(n,c/n)$ remains asymptotically the
same for all values of the constant $c>0$ and appears thus to be quite
insensitive to the growth of $c>0$. This stays in a sharp contrast with
many other properties of random graphs such as the appearance of the
giant component (all components of $G(n,c/n)$ are a.s. at most
logarithmic in size for $c<1$, while for $c>1$ $G(n,p)$ contains a.s.
one component of a linear size and the rest are $O(\log n)$) or
planarity ($G(n,c/n)$ is a.s. planar for $c<1$ and a.s. non-planar for
$c>1$). 

Another related problem is to investigate the spectrum of the Laplacian
of a random graph $G(n,p)$. For a graph $G$, the Laplacian $L=L(G)$ is
defined as $L=D-A$, where $A$ is the adjacency matrix of $G$ and $D$ is
the diagonal matrix whose diagonal entries are degrees of corresponding
vertices. For any graph $G$, the Laplacian $L(G)$ is easily seen to be a
real symmetric matrix with non-negative eigenvalues, the smallest of
them being zero. One may study the so-called spectral gap (the smallest
positive eigenvalue of the Laplacian) of random graphs
$G(n,p)$ for various values of $p(n)$.

The methods of this paper can be possibly applied to the study of the
spectrum of dilute random matrices. A dilute random matrix $A$
is defined by 
\begin{eqnarray*}
A_{i,j}&=&a_{i,j}b_{i,j}, \quad 1\le i\le j\le n\\
A_{j,i}&=&A_{i,j}, \quad 1\le i<j\le n .
\end{eqnarray*}
where $a_{i,j}$ are jointly independent random variables with zero mean
and variance 1, and $b_{i,j}$ are also jointly independent and
independent from $\{a_{i,j}\}$, where $b_{i,j}=1$ with probability
$p=p(n)$ and $b_{i,j}=0$ with probability $1-p$. In other words, the
dilute random matrix is obtained by replacing each entry of a matrix
from the so-called Wigner ensemble by zero independently with
probability $q=1-p$. As such, it unifies the notions of the Wigner
random matrices and random graphs. Khorunzhy proved in \cite{Kho01a} that
the spectral norm of the dilute random matrix is asymptotically equal to
$2\sqrt{np}$ in the case $p(n)\gg \log n/n$ (under some additional
technical assumptions) and it asymptotically much larger than
$\sqrt{np}$ for $p(n)\ll \log n$. It would be quite interesting to
determine the asymptotic behavior of the spectral radius of the dilute
random matrix for the case of small values of $p(n)$.


\begin{thebibliography}{99}
\bibitem{Alo98} N. Alon,  Spectral techniques in graph algorithms,
in: {\em LATIN'98: theoretical informatics} (Campinas, 1998), 206--215, 
Lecture Notes in Comput. Sci., Vol. 1380, Springer, Berlin, pp.
206--215.
\bibitem{AS} N. Alon and J. H. Spencer, {\bf The probabilistic
method}, Wiley, New York, 1992.
\bibitem{B} B. Bollob\'as, {\bf Random graphs}, Academic Press,
New York, 1985.
\bibitem{Chu97} F. R. K. Chung, {\bf Spectral graph theory}, Regional
Conference Series in Mathematics 92, Amer. Math. Soc., Providence, 1997.
\bibitem{CveDooSac95} D. Cvetkovi\'c, M. Doob and H. Sachs,  
{\bf Spectra of graphs. Theory and 
applications}, $3^{rd}$ ed., Johann Ambrosius Barth, Heidelberg, 1995.
\bibitem{FurKom81} Z. F\"uredi and J. Koml\'os, {\em The eigenvalues of
random symmetric matricesi}, { Combinatorica} 1 (3), 233--241 (1981).
\bibitem{JanLucRuc00} S. Janson, T. \L uczak and A. Ruci\'nski, {\bf
Random graphs}, Wiley, New York, 2000.
\bibitem{Kho01a} A. Khorunzhy, {\em Sparse random matrices: spectral
edge and statistics of rooted trees}, Adv. Appl. Probab. 33 (2001),
1--18.
\bibitem{Kho01b} A. Khorunzhy, {\em Stochastic version of the Erd\H
os-R\'enyi limit theorem}, Preprint {\em xxx.lanl.gov} :
math.PR/0103043. 
\bibitem{KhoVen00} A. Khorunzhy and V. Vengerovsky, {\em On asymptotic
solvability of random graph's Laplacians}, Preprint {\em xxx.lanl.gov} :
math-ph/0009028.
\bibitem{KriVu00} M. Krivelevich and V. H. Vu, {\em On the concentration
of eigenvalues of random symmetric matrices}, Preprint
{\em xxx.lanl.gov} : math-ph/0009032.
\bibitem{L} L. Lov\'asz,
{\bf Combinatorial Problems and Exercises},
North-Holland, Amsterdam, 1993.
\bibitem{Meh91} M. L. Mehta, {\bf Random matrices}, Academic Press,
New York, 1991.

\end{thebibliography}
\end{document}